\input amstex 
\documentstyle{amsppt}
\input bull-ppt
\keyedby{bull273/PAZ}

\topmatter
\cvol{26}
\cvolyear{1992}
\cmonth{April}
\cyear{1992}
\cvolno{2}
\cpgs{280-287}
\title A period mapping in universal Teichm\"uller 
space\endtitle
\author Subhashis Nag\endauthor
\address The Institute of Mathematical Sciences,
C.I.T. Campus,
Madras 600 113, INDIA\endaddress
\ml nag\@imsc.ernet.in\endml
\date May 21, 1991 and, in revised form, August 20, 
1991\enddate
\thanks This paper was presented as an invited lecture on 
May 27, 1991, 
at the second International Symposium on Topological and 
Geometrical
Methods in Field Theory, Turku, Finland, (May 27--June 1, 
1991). Sponsored
by Academy of Finland et al\endthanks
\abstract In previous work it had been shown
that the remarkable homogeneous space $M=
\operatorname{Diff}(S^1)/\operatorname{PSL}
(2,\Bbb{R})$ sits as a complex analytic
and K\"ahler submanifold of the Universal Teichm\"uller
Space. There is a natural immersion $\Pi$ of $M$
into the infinite-dimensional version (due to Segal)
of the Siegel space of period matrices. That
map $\Pi$ is proved to be injective, equivariant,
holomorphic, and K\"ahler-isometric (with respect
to the canonical metrics). Regarding a period
mapping as a map describing the variation of complex
structure, we explain why $\Pi$ is an infinite-dimensional
period mapping.\endabstract
\endtopmatter

\document
\heading {}  Introduction\endheading
\par
Let $\operatorname{Diff}(S^1)$
denote, as usual, the infinite-dimensional Lie
group of orientation preserving $C^\infty$ diffeomorphisms
of the unit circle. In a previous paper \cite{12} we showed
that the canonical complex analytic structure and K\"ahler
metric, $g$, that exist on the homogeneous space 
$M=\operatorname{Diff}
(S^1)/\operatorname{\text{M\"ob}}(S^1)$ {\it coincide 
exactly\/}
with the canonical Ahlfors-Bers complex structure and the
(generalised) Weil-Petersson K\"ahler structure of the
Universal Teichm\"uller Space $T(1)$. This identification
was via the natural injection {\rm(I)} of $M$ into $T(1)$,
where $T(1)$ is thought of as the homogeneous space of all
quasi-symmetric homeomorphisms of $S^1$ modulo the 
M\"obius subgroup.
\par
In this paper we study another natural
embedding $(\Pi)$ of $M$ into the
infinite-dimensional version, $D_\infty$,
of the Siegel disc \RM(or
generalised upper halfspace\RM) of period matrices. These 
matrices
(suitable symmetric Hilbert-Schmidt operators) comprising 
$D_\infty$ are
symmetric with complex entries, so $D_\infty$ is a complex 
manifold
(infinite-dimensional bounded domain). As in finite 
dimensions, so
also here the Siegel symplectic (K\"ahler) metric, $h$, 
exists on this
homogeneous manifold. One main purpose here is to announce 
that,
like \RM I, the map $\Pi$ is completely natural \RM(in 
that it
respects all the structures\RM)\RM; namely, $\Pi$ {\it is 
an equivariant,
holomorphic, and K\"ahler isometric immersion of $M$ into 
$D_\infty$}.
Thus, combining this with the results proved by this 
author in
\cite{12, Part II}, one can assert that on the complex 
submanifold $M$
of the Teichm\"uller space, the Weil-Petersson metric 
coincides with the
Siegel symplectic metric.
\par
Considering $M$ as the submanifold of {\rm``}smooth{\rm''} 
complex
structures in the Teichm\"uller space, we can interpret 
the mapping
$\Pi$ as a period mapping that associates to each point of 
$M$ its
corresponding matrix of periods. The idea is to look at 
the period
map in P. Griffiths's way as a map describing the 
variation of Hodge
and complex structure. We will provide an exposition.
\par
As explained in \cite{12}, that work was intimately 
related to
and motivated by string physics. Once again, physicists 
have been
looking at the map $\Pi$ and have made several claims. The 
various
infinite-dimensional spaces appearing in the set-up are all
somehow {\it universal moduli spaces\/} and may turn out 
to be
important in a nonperturbative formulation of string theory.
See Nag \cite{10}.
\heading 1. The map $\Pi$\endheading
\par
The method of defining $\Pi$ stems
from a faithful representation of
$\operatorname{Diff}(S^1)$ in the
infinite-dimensional symplectic group,
say $\operatorname{Sp}$, that was
introduced by G. Segal \cite{13, \S5}.
(Segal's motivation is to obtain
the metaplectic representation of
Sp, but that need not concern us here.)
Recall that if a (real) vector space $V$
is equipped with a skew-symmetric, nondegenerate
bilinear form, $S$, then the linear automorphisms
$T\in\operatorname{GL}(V)$ that preserve
the form $S$ (viz., $S(Tv_1,Tv_2)=S(v_1,v_2)$
for all $v_1$, $v_2$ in $V)$ constitute the symplectic
subgroup $\operatorname{Sp}(V,S)\subset \operatorname{GL}
(V)$.
\par
Consider the vector space
$$
V=C^\infty \operatorname{Maps}
(S^1,\Bbb{R})/\Bbb{R}\text{\ (constant maps)}.
\tag1
$$
The quotienting by the one-dimensional subspace of 
constant maps
means that one is normalizing the 0-th Fourier coefficient 
to vanish.
$V$ is equipped with a natural nondegenerate 
skew-symmetric bilinear
form $S\: V\times V\rightarrow\Bbb{R}$ (Segal \cite{13}), 
given by
$$
S(\sigma,\tau)=\int_{S^1} 
\sigma(\theta)(d\tau/d\theta)\,d\theta.\tag2
$$
Then $\operatorname{Diff}(S^1)$ acts on $V$,
by the obvious precomposition, preserving this
bilinear form $S$. {\it Hence the diffeomorphism group 
becomes
a subgroup of the real symplectic automorphisms of $(V,S)$}.
Actually one deals with a restricted subgroup of the 
symplectic
automorphisms, denoted $\operatorname{Sp}_0(V)$, which 
already contains
$\operatorname{Diff}(S^1)$ (see \cite{13}). (This is an 
analytic detail
to simplify the infinite-dimensional considerations, since
$\operatorname{Sp}_0(V)$ elements have a natural norm.)
\par
To proceed further we recall the definition of a {\it 
positive 
polarization\/} of the space $V$ with respect to $S$. Let $V
_{\Bbb{C}}$ denote the complexification of $V$. Then a 
decomposition
$V_{\Bbb{C}}=W\oplus\overline W$ such that (the 
complexification of)
$S$ takes zero values on arbitrary pairs from $W$ is called
a polarization and $W$ is called a (maximal) {\it 
isotropic subspace\/} for
$S$. The assignment
$$
\langle w_1,w_2\rangle=iS(\bar w_1,w_2)\tag3
$$
is a Hermitian inner product on $W$, and the decomposition 
is a
{\it positive polarization\/} if this is positive 
definite. In that case
$W$, its conjugate, and hence $V_{\Bbb{C}}$ itself, can be 
completed to
Hilbert spaces with respect to this Hermitian inner 
product. We will henceforth
identify a positive polarization with the isotropic subspace
$W$ determining it. Here is the {\it fiducial\/} positive 
polarization
for $(V,S)$: take $W=W_+$ to be the subspace of 
$V_{\Bbb{C}}=C^\infty
\operatorname{Maps}(S^1,\Bbb{C})/\Bbb{C}$ consisting of 
those mappings having
only positive index Fourier components. Note that 
$\overline W_+
=W_-$, and also the fundamental fact that the image under 
\RM(the
$\Bbb{C}$-linear extension of\RM) a symplectic 
automorphism of a positive
isotropic subspace is again such a subspace. In fact, 
$\operatorname{Sp}
(V,S)$ acts transitively on the space of all positive 
polarizations, and the
stabilizer subgroup at $W$ is evidently identifiable with 
the unitary group
$U(W,\langle\ ,\ \rangle)$. It follows that the 
homogeneous space
$\operatorname{Sp}/U$ can be identified with the family of 
positive
polarizations of $V$. Passing to the restricted subgroup 
$\operatorname{Sp}_0$
one finds that $\operatorname{Sp}_0/U$ represents a class 
of (``bounded'')
positive polarizations of $V$, say 
$\operatorname{Pol}_0(V)$, and either
of these spaces can be easily identified with the Siegel
disc of infinite dimension (Segal \cite{13, p.\ 316}),
$$
\aligned
D_\infty=&\{\text{All Hilbert-Schmidt complex linear 
operators $Z\:
W_+\rightarrow W_-$}\\
&\hphantom{\{}\text{such that $Z$ is symmetric (w.r.t.\ $S)$
and $I-Z\overline Z$ is positive definite\}.}\endaligned
\tag4
$$
For instance, the identification between $D_\infty$ and
$\operatorname{Pol}_0(V)$ is by associating to $Z\in 
D_\infty$
the positive isotropic subspace $W$ that is the graph of the
operator $Z$. (The origin in $D_\infty$ corresponds to the
fiducial polarization $W_+$.) In finite dimension (genus $g$
Riemann surfaces) this identification corresponds to 
taking $W$
to be the row span of the ``full period matrix'' $(Z,I)$; 
then $W$
is a positive isotropic $g$-dimensional complex subspace of 
$\Bbb{C}^{2g}$ (equipped with standard skew form). The 
current
set-up thus agrees well with the classical case, and we 
have in
infinite dimension the Siegel disc appearing in the various
incarnations
$$
D_\infty=\operatorname{Sp}_0/U=\operatorname{Pol}_0(V)%
\hookrightarrow
\operatorname{Gr}(W_+,V_{\Bbb{C}}).\tag5
$$
Here the subspaces $W\in\operatorname{Pol}_0(V)$ comprise 
a complex
analytic submanifold of the complex Grassmannian of 
subspaces
of $V_{\Bbb{C}}$ that are graphs of (not necessarily 
symmetric)
linear operators of Hilbert-Schmidt type from $W_+$ to 
$W_-$. This
complex structure on $\operatorname{Pol}_0(V)$ (identified 
as explained
with $D_\infty)$ is easily checked to be the same as the 
complex structure
of $D_\infty$ as a family of (symmetric with respect to 
$S)$ complex
operators.
\par
We are finally ready to define $\Pi$.
\thm{Lemma}
In Segal\/\RM's representation 
$\operatorname{Diff}(S^1)\rightarrow
\operatorname{Sp}_0(V)$, the diffeomorphisms that map into 
the
unitary subgroup $U$ of $\operatorname{Sp}_0$ comprise 
precisely
the \RM3-parameter subgroup of M\"obius transformations
on $S^1$. Consequently one obtains an injective 
equivariant mapping
$$
\Pi\:\operatorname{Diff}(S^1)/\text{\rm M\"ob}(S^1)
\rightarrow \operatorname{Sp}_0/U=D_\infty.
\tag6
$$
\ethm
\par
The proof boils down to showing that a holomorphic self-map
of the disc that extends continuously to the boundary as a
diffeomorphism of $S^1$ must be a M\"obius transformation. 
This
follows, for example, by the argument principle. That 
$\Pi$ is
equivariant (with respect to corresponding right 
translations
on  domain and target) is easily verified.
\heading 2. Holomorphy of $\Pi$\endheading
\thm{Theorem}
$\Pi$ is a holomorphic immersion.
\ethm
\par
The method is to compute the derivative of $\Pi$ at the 
origin of
$M$ and show that it is a complex linear map of the 
tangent space of
$M$ to the tangent space of the Grassmannian above. By 
equivariance
it is enough to verify this at the origin.
\par
Let me first explain the canonical complex structure on 
the homogeneous
space $M$. The Lie algebra of the Lie group 
$\operatorname{Diff}(S^1)$
is the algebra of smooth $(C^\infty)$ real vector fields 
on $S^1$.
The tangent space at the origin (i.e., coset of the 
identity)
of $M$ corresponds to those smooth vector fields 
$v=u(\theta)
\partial/\partial \theta$ such that in the Fourier expansion
$$
u(\theta)=\sum^\infty_{m=-\infty} u_m e^{im\theta},
\tag7
$$
one has $u_{-1}=u_0=u_1=0$. Note that reality of $v$ implies
$\overline{u_m}=u_{-m}$. One easily sees that the three 
middle Fourier
coefficients vanish because an infinitesimal M\"obius 
transformation allows
one to normalize precisely these three coefficients.
\par
The canonical homogeneous almost-complex structure $J$ on 
$M$ is now
defined simply by {\rm``}conjugation of Fourier 
Series,{\rm''} namely,
$$
J\left(u(\theta)\frac{\partial}{\partial\theta}\right)
=u^\ast(\theta) \frac{\partial}{\partial\theta},
\tag8
$$
where
$$
u^\ast(\theta)=-\sum^\infty_{m=-\infty}
i\operatorname{sgn}(m)u_m e^{im\theta}.
\tag9
$$
This defines $J$ at the origin on the homogeneous 
Fr\'echet manifold
$M$; elsewhere $J$ is transported by right translations. 
It is easy
to check (at least at the formal level), that this almost 
complex
structure on $M$ is integrable. Thus $M$ is canonically a 
homogeneous
complex manifold. The complex structure above arises 
naturally
from various points of view, for instance, from the 
Kirillov-Kostant
coadjoint orbit theory and also in work of Pressley et al. 
See
references \cite{1--4, 16, 17}.
\par
Now, the tangent space to the Grassmannian at a point $H$ 
can
be canonically identified with the linear space 
$\operatorname{Hom}
(H,V_{\Bbb{C}}/H)$ (only those homomorphisms satisfying a 
Hilbert-Schmidt
boundedness condition are included). We now exhibit the 
derivative of $\Pi$ from the
tangent space of $M$ to the tangent space of the 
Grassmannian.
\thm{Proposition}
The derivative of $\Pi$ at the origin of $M$
carries the vector field 
$v=u(\theta)\partial/\partial\theta$
to the following homomorphism $\eta_v\in\operatorname{Hom}
(W_+,W_-)$\RM:
$$
\eta_v(\sigma)=\pi_-(u(\theta)\sigma'(\theta))\quad
\text{for any }\sigma\in W_+,\tag10
$$
where $\pi_-$ denotes the projection of $V_{\Bbb{C}}$ onto 
$W_-$
and prime denotes differentiation with respect to $\theta$.
\ethm
\par
The proof of (10) is straightforward after one has properly
identified the various spaces involved.
\par
To establish the holomorphy of $\Pi$ one needs to prove that
{\it the homomorphism $\eta_{Jv}$ is $i$ times $\eta_v$.}
Recall that the negative index Fourier coefficients of 
$Jv$ are just the
corresponding coefficients of $v$ multiplied by $i$. On 
the other hand,
any $\sigma$ from $W_+$ (therefore any $\sigma')$ has only 
positive
index Fourier components. Application of formula (10), 
therefore,
gives what we wish.
\heading 3. $\Pi$ is an isometry\endheading
\par
There is a unique homogeneous K\"ahler metric, $g$, 
available on
$M$. In Part II of \cite{12} I had proved that this metric 
is precisely
the Weil-Petersson metric on the Teichm\"uller spaces.
\par
Here is one way to explain how $g$ arises. Since $M$ is a
coadjoint orbit in the coadjoint representation of the 
Virasoro
group, (see Witten \cite{16}), one has, as always, the 
natural
Kirillov-Kostant symplectic form on these homogeneous 
manifolds.
This symplectic form is compatible with the complex 
structure seen
above, and we thus get the K\"ahler structure $g$ on $M$. 
In a less
sophisticated fashion, it is actually easy to see that $g$ 
is the
unique homogeneous K\"ahler metric (up to a numerical 
scaling
parameter) possible on $M$. The number 26, namely, the 
dimension of
bosonic spacetime, appears entrenched in the curvature of 
$g$. See
calculations of Bowick, Rajeev, Lahiri, and Zumino in the 
references.
Explicitly, this canonical $g$, at the origin of $M$ has 
the Hermitian
pairing (see \cite{12, Part II})
$$
g(v,w)=\left[\sum^\infty_{m=2} v_m 
\overline{w_m}(m^3-m)\right].\tag11
$$
Here $v$ and $w$ are two smooth real vector fields on 
$S^1$ representing
tangent vectors to $M$ (see equation (7)), and the $v_m$, 
$w_m$ are the
respective Fourier coefficients. Elsewhere on $M$, $g$ is 
of course
transported by right translations isometries.
\par
The {\it Siegel symplectic metric\/}, $h$, (note Siegel 
\cite{15}) exists on
each finite-dimensional Siegel disc $D_g$ as the unique 
K\"ahler metric
on that bounded domain for which the full holomorphic 
automorphism group
$\operatorname{Sp}(2g,\Bbb{R})$ acts as isometries. It 
generalises without
trouble to Segal's infinite-dimensional version $D_\infty$ 
explained above.
The pairing given by $h$ on the symmetric elements of 
$\operatorname{Hom}
(W_+,W_-)$, i.e., on the tangent space at the origin of 
$D_\infty$, turns
out, after some work, to be
$$
h(\phi,\psi)=\operatorname{trace\ of} \phi\circ\bar 
\psi,\tag12
$$
where the conjugate of $\psi$ maps $W_-$ to $W_+$, so 
taking trace makes
sense (and gives a finite answer).
\thm{Theorem}
$\Pi$ is an isometric immersion of $(M,g)$ into 
$(D_\infty,h)$.
\ethm
\par
Again utilising equivariance and the fact that $g$ and $h$ 
are
homogeneous metrics we only need check at the origin. One
calculates the trace of $\eta_v$ composed with conjugate 
of $\eta_w$
and finds interestingly that the answer is as predicted by 
formula
(11). The computation is finally a finite-series 
summation. See \cite{10}
for details.
\heading 4. Dependence of $\Pi$ on Beltrami 
coefficients\endheading
\par
Here is the explicit form of the period matrix in $D_\infty$
associated to a given Beltrami coefficient $\mu$ 
(representing a Teichm\"uller
point $[\mu])$. The symplectic automorphism $T_\mu$ of $V$ 
arising from
$\mu$ is $h\mapsto h\circ w_\mu$, $h\in V$. Here $w_\mu$ 
{\it is the
$\mu$-conformal self-homeomorphism of the unit disc 
$\Delta$}, and we are
using its boundary homeomorphism. The matrix for $T_\mu$ 
(complexified)
in the standard orthonormal bases $\{e_k=e^{ik\theta} 
/\sqrt{|k|}\}$
for $V_{\Bbb{C}}=W_+\oplus W_-$ has the form
$$
T_\mu=\bmatrix A&B\\
\overline B& \overline A\endbmatrix.\tag13
$$
Here $A\: W_+\rightarrow W_+$, $B\: W_-\rightarrow W_+$,
$A=((a_{pq}))$, $B=((b_{pq}))$.
\par
We obtain the fundamental formulas
$$
\align
a_{pq}&=\frac{1}{2\pi} \frac{\sqrt{p}}{\sqrt{q}} 
\int^{2\pi}_0
 (w_\mu(e^{i\theta}))^q e^{-ip\theta}\,d\theta,\qquad p,q>0;
\tag14\\
b_{rs}&=\frac{1}{2\pi}\frac{\sqrt{r}}{\sqrt{s}}\int^{2\pi}_0
 (w_\mu(e^{i\theta}))^{-s} e^{-ir\theta}\,d\theta,\qquad 
r,s>0.
\tag15
\endalign
$$
Recall that the symplectic group acts by holomorphic 
automorphisms
on $D_\infty\: Z\mapsto (\overline A\circ Z+\overline 
B)\circ
(B\circ Z+A)^{-1}$ where the symplectic matrix ``$T$'' is 
in the block
form (13). The zero matrix in $D_\infty$ corresponds to 
$\mu=0$,
and we see that
$$
\Pi([\mu])=\overline B\circ A^{-1}\tag16
$$
with $A$ and $B$ given as above.
\subheading{First variation of the period matrix}
That the matrix (16) varies holomorphically
with $\mu$ is already nontrivial (remember that $w_\mu$
varies only real analytically with $\mu$). Since
we know the perturbation formula for $w_\mu$, we
can actually write down the infinitesimal variation of
the matrix entries of $\Pi(\mu)$ as $\mu$ varies.
\thm{Proposition}
For $\mu\in L_\infty(\Delta)$ one has, as $t\rightarrow0$,
$$
\Pi([t\mu])_{rs}=t\cdot\pi^{-1}
\cdot\sqrt{rs}\iint_\Delta
\mu(z) z^{r+s-2}\,dx\,dy+O(t^2).\tag17
$$
\ethm
\rem{Remark}
The proposition corresponds precisely to the well-known
{\it Rauch variational formulas} for finite genus period 
mappings. In
fact, the holomorphic 1-forms on $\Delta$ are generated by 
$\phi_r
=z^{r-1}$, $r=1,2,\dotsc$; so (17) is in harmony with the 
Rauch
formula ``\<$\int\!\int\mu\phi_r\phi_s$'' for the 
variation of the $rs$
entry. (See \cite{11, Chapter 4}.)
\endrem
\heading 5. $\Pi$ as a period mapping\endheading
\par
The map $\Pi$ links up with the variation of Hodge
structure following P. Griffiths \cite{5, 6}. The idea is to
interpret the vector space $V$ as the ``first cohomology 
(with
real coefficients) of a surface of infinite genus.'' The 
skew
form $S$ can then be recognised as the cup product map on 
this
cohomology vector space.
\par
Now, a ``complex structure'', that is a point from the
universal Teichm\"uller space, will produce a 
decomposition of the
complexification of $V$ into its $H^{1,0}$ and $H^{0,1}$ 
subspaces.
The $H^{1,0}$ subspace is a positive isotropic subspace 
for the cup
product skew form because the condition of being positive 
isotropic is an
implementation of the bilinear relations of Riemann. The 
decomposition so
obtained is indeed a positive polarization of $V$. For the 
base complex
structure ``0'', we take $H^{1,0}$ to be $W_+$. Thus $\Pi$
becomes realised as an association between the smooth 
complex
structures and their corresponding $H^{1,0}$ subspaces, as 
is
expected of a period mapping in Griffiths's set-up. That 
$\Pi$
is injective (as we saw) is then {\it Torelli\RM's 
Theorem\/} in this
context.
\par
One can also prove that $\Pi$ keeps track of the varying 
space of
holomorphic 1-forms as the complex structure varies. That is
precisely the characteristic of the classical period 
mappings. See
\cite{10} for this point of view.
\subheading{The infinitesimal Schottky locus}
It is clear that the image $\Pi(M)$
in $D_\infty$ is the analog of the
Schottky locus of Jacobians in the Siegel
disc $D_g$. Using the formula (10) for $d\Pi$
we get a pretty description of the tangent to
the ``Schottky locus'' $\Pi(M)$:
\thm{Theorem}
The tangent space to $D_\infty$ at the origin
is canonically identifiable with the complex symmetric
matrices $((\lambda_{pq}))$ $(p,q=1,2,3,\dotsc)$
satisfying the Hilbert-Schmidt condition
$\sum^\infty_{p=1} \sum^\infty_{q=1} 
|\lambda_{pq}|^2<\infty$.
The image of $d\Pi$ at the origin consists of those
that are of the form
$$
\lambda_{pq}=i\sqrt{pq} a_{(p+q)},
\tag18
$$
where the sequence $(a_2,a_3,a_4,\dotsc)$ is an arbitrary
complex sequence whose $n$\<th term goes to zero faster than
$|n|^{-k}$ for any $k>0$.
\ethm
\rem{Remark}
In finite genus $g$,
the $(3g-3)$-dimensional Teichm\"uller space
sits via $\pi_g$ in the $g(g+1)/2$-dimensional
$D_g$. Here also we see that the number
of independent entries in a $D_\infty$
tangent vector is growing quadratically with the
block size, whereas the ones arising from
Teichm\"uller points are determined just by their
first column, i.e., grows linearly with block
size. In this aspect also $\Pi$ generalises the
$\pi_g$.
\endrem
\rem{Remark}
Since $D_\infty$ sits within the infinite-dimensional
Grassmannian, the theory (of the tau function) developed
in \cite{14} by Segal-Wilson suggests that it may be
possible to describe the image of $\Pi$ by solutions of
$K-P$ equations; namely, an infinite-dimensional form of
the Novikov conjecture could be valid.
\endrem
\rem{Remark}
In Part II of \cite{12} the present author had shown a
form of Mostow rigidity that implied that the
Teichm\"uller spaces $T_g$ of compact Riemann surfaces
of genus $g$ sit within the universal $T(1)$ cutting
transversely the submanifold $M$. That result disallows
the possibility of relating the finite-dimensional period
maps $\pi_g\: T_g\rightarrow D_g$ to $\Pi$ by simple 
restriction
of domains. Nevertheless, it was possible in Part II of 
\cite{12} to relate the Weil-Petersson metric on $T_g$
to the canonical metric $g$ of $M$ by a regulation of
improper integrals. Some similar analysis may relate $\pi_g$
to $\Pi$.
\endrem
\par
In any case the result of \S3 raises the following 
interesting question:
How is $\pi^\ast_g$(Siegel symplectic metric) related to 
the Weil-Petersson
metric on $T_g$? Certainly they do not coincide in low 
genus! A connection
between the volume form induced on $T_g$ by the Siegel 
metric and volume
forms arising from the complex geometry of $T_g$ was shown 
in Nag
\cite9.
\par
Hong-Rajeev \cite7 have considered the mapping $\Pi$ as a 
map from
the space of univalent functions (Riemann maps) to 
$D_\infty$; they have
made several claims but have offered no proofs. They 
appear to be
thinking of points of $\operatorname{Sp}_0/U$ as complex 
structures on $V$
rather than as polarizations. It seems to us that this 
does not quite fit
requirements. Also, if one wishes to look at $\Pi$ as a 
map on
univalent functions, namely, as a map on domains bounded 
by smooth Jordan
curves, then one is forced to use the {\it conformal 
welding 
correspondence\/} to relate to $\operatorname{Diff}(S^1)/
\operatorname{\text{M\"ob}}(S^1)$. (See \cite{12, Part II} 
and Katznelson-Nag-Sullivan
\cite8). Using the welding it is possible for us to write 
formulas for
$\Pi$ as a map on coefficients of Riemann maps (which are 
essentially nothing
but Bers embedding holomorphic coordinates); however, we 
are unable to see
how the formulas presented by Hong and Rajeev can be valid.
\par
More details on this work will appear in \cite{10} and
elsewhere.

\enddocument